\documentclass[11pt,a4paper,leqno]{article}
\usepackage{latexsym}
\usepackage{amsfonts}
\usepackage {amsbsy}
\usepackage {amsmath}
\usepackage{amssymb}
\newtheorem{defn}{Definition}[section] 
\newtheorem{thm}[defn]{Theorem} 
\newtheorem{prop}[defn]{Proposition} 
\newtheorem{lem}[defn]{Lemma} 
\newtheorem{cor}[defn]{Corollary}

\newcommand \psh {plurisubharmonic }
\newcommand \proof { Proof: }
\newcommand \C{\mathbb C}
\newcommand \N{\mathbb N}

\newcommand \R{\mathbb R}
\newcommand \B {\mathbb B}
\newcommand \fin{$\blacktriangleright $\\}
\newcommand \vep  {\varepsilon}
\newcommand  \vphi {\varphi}
\newcommand \mcal  {\mathcal}

\newcommand \W {\Omega}
\newcommand \E {\mathcal E}
\newcommand \F {\mathcal F}
\newcommand \w {\omega}

\newcommand \x {\times}
\newcommand \Sub {\Subset}
\newcommand \sub{\subset}
\newcommand \ove{\overline}
\newcommand \sm{\setminus}

\newcommand \Om{\Omega }

\newcommand \MA{Monge-Amp\`ere }

\newcommand \fii{\varphi }

\newcommand \ep{\epsilon }

\newcommand \al{\alpha }

\newcommand \du{(dd^c u )^n }

\newcommand \dduj{(dd^c u_j )^n  }
\newcommand \ddvj{(dd^c v_j )^n  }
\newcommand \ddwj{(dd^c w_j )^n  }

\newcommand \lp{\mathcal  L }

\newcommand \ffa{\mathcal  F ^a}
\numberwithin{equation}{section}
 \title{Subextension of  \psh  functions   \\
 with weak singularities}

 \author{U. CEGRELL{\footnote {partially supported by the Swedish Research
  Council contact no 621-2002-5308}}, \
 S. KOLODZIEJ{\footnote{partially supported by KBN grant 1 P03A 03727}} \ and A. ZERIAHI} 
 \date{}
\begin{document}
\maketitle
\section{Introduction} 
E. Bedford and D. Burns ([Be-Bu]) and later U. Cegrell ([Ce 1]) 
proved around 1978  that any smooth bounded domain satisfying 
certain boundary conditions is a domain of existence of a 
plurisubharmonic function.

However since \psh functions occur in complex analysis 
through inequalities, it is more natural to ask  for the 
subextension problem.

El Mir gave in 1980 an example of a \psh function on the 
unit bidisc for which the restriction to any 
smaller bidisc admits no subextension to the 
whole space (see [El]). He also proved that, after 
attenuatting the singularities of a given \psh 
function by composition with a suitable convex increasing 
function, it is possible to obtain a global 
subextension.

Later Alexander and Taylor gave in 1984 a generalization of this 
result with a more effective and simple proof (see [Al-Ta]).

On the other hand, Fornaess and Sibony pointed
out in 1987 that for a ring domain in $\C^{2},$ there 
exists a \psh function which admits no 
subextension inside the hole (see [Fo-Sib]).  

 Finally E. Bedford and B.A. Taylor proved in 1988
 that any smoothly bounded domain in $\C^n$ is a domain of existence of a 
smooth \psh function (see [Be-Ta 3]).
 
Recently, the first and the last authors proved that \psh functions with 
uniformly bounded Monge-Amp\`ere mass on 
a bounded hyperconvex domain always admit a \psh 
subextension to any larger hyperconvex domain (see [Ce-Ze]). 
 
 Here we want to prove several results showing that  \psh functions with
 various bounds on their Monge-Amp\`ere masses on 
a bounded hyperconvex domain always  admit global \psh 
subextension with logarithmic growth at infinity.

Let us describe more precisely the content of the article.

In section 2 we recall basic definitions concerning the Cegrell class $\mcal E (\W)$ of
\psh functions of locally uniformly bounded Monge-Amp\`ere masses on a hyperconvex domain
$\W \Sub \C^n$ and its subclasses of \psh functions of finte energy. Then we give a characterization in
terms of some capacity of
functions from the class $\mcal E (\W)$  which provides several examples of functions in this class.

In section 3, we give almost sharp estimates on the size of sublevel sets of \psh functions
in various subclasses of the class $\mcal E (\W).$

In section 4, we give a generalization of
Alexander-Taylor's subextension  theorem which implies, using results from
section 3,
 that \psh  functions of finite energy 
in the sense of Cegrell admit a global subextention with logarithmic 
growth of arbitrary small logarithmic type. 

 Finally in section 5, using recent results from the theory of Monge-Amp\`ere equation on compact
 K\"ahler manifolds due to the second author, we prove two results about global
 subextension of \psh functions of uniformly bounded masses on a hyperconvex
 domain by \psh function with logarithmic growth on $\C^n$ with a
 well defined global Monge-Amp\`ere measure in some
 cases.

 A part of this research was done during the visit of the second author in
Laboratoire Emile Picard of the University Paul Sabatier (Toulouse) in
 january 2004. He would
like to thank the staff of the Laboratoire for its hospitality and stimulating discussions.
 
\section {Plurisubharmonic functions with locally uniformly bounded Monge-Amp\`ere masses}  
Let us first recall some definitions from ([Ce 2], [Ce 3]). Let 
$\W \Sub  \C^n$ be a hyperconvex domain. We denote 
by $\mcal E_{0} (\W)$ the set of negative and bounded \psh 
functions $ \vphi$ on $\W$ which tend to zero at the boundary 
and satisfy $\int_{\W} (dd^c \vphi)^n < + \infty.$
Then for each $p \geq 1$ define  $\mcal E_{p} (\W)$ to be the 
class of \psh functions $\vphi$ on $\W$ such that there exists 
 a decreasing sequence of \psh functions $(\vphi_{j})$ from the 
 class $\mcal E_{0} (\W)$ which converges to $\vphi$ such that 
 \begin{equation}
        \sup_{j} \int_{\W} (-\vphi_{j})^p 
 (dd^c \vphi_{j})^n < + \infty.
        \label{eq:Ep}
 \end{equation}
If the sequence $(\vphi_{j})$ can be choosen 
so that $ \sup_{j}  \int_{\W}(dd^c \vphi_{j})^n 
< + \infty,$ then we say that 
 $\vphi \in \mcal F_{p} (\W).$

 Let us denote by $\mcal F (\W)$ the set of all 
$\vphi \in PSH (\W)$ such that there exists a 
sequence  $(\vphi_{j})$ of \psh function 
in $\mcal E_{0} (\W)$ such that 
$\vphi_{j} \searrow \vphi$ and $\sup_{j} \int_{\W} (dd^c 
\vphi_{j})^n < + \infty.$ We also need the subclass
$\mcal F^a (\W)$ of functions from $\mcal F (\W)$ whose Monge-Amp\`ere
measures put no mass on pluripolar subsets of $\W.$ 
 By Cegrell [Ce 2] we have $\mcal E_{0} (\W) 
\sub \mcal F_{p} (\W) \sub \mcal E_{p} (\W) \sub \mcal F^a (\W) \sub \mcal F (\W),$ 
for any $p \geq 1.$
 It is also proved in ([Ce 2], [Ce 3]) that the 
complex Monge-Amp\`ere operator is well 
defined for a function $\vphi \in \mcal F (\W)$ as the weak*-limit of the sequence of 
measures $(dd^c \vphi_{j})^n,$ where 
$(\vphi_{j})$ is any decreasing sequence of 
\psh functions from the class $\mcal E_{0} (\W)$ 
which converges to $\vphi$ and satisfying the required conditions in the definition.

  Finally we denote by $\mcal E (\W)$ the set of \psh functions 
  which are locally in $\mcal F (\W).$ 
Then the complex Monge-Amp\`ere operator is well defined on the class $\mcal E
  (\W)$ (see [Ce 3]).\\
  The following result will be useful for later refrence; it has been proved in [Ce
  3] for $p \geq 1.$
 \begin {lem} Let $\vphi \in \mcal E_{p} (\W),$ if $p > 0$ 
and $\vphi \in \mcal F (\W)$ if $p = 0.$ 
Then the pluricomplex $p-$energy of $\vphi$ 
defined by the following formula
\begin{equation}
 e_{p} (\vphi) :=  \int_{\W} (-\vphi)^p 
 (dd^c \vphi)^n
        \label{}
\end{equation}
is finite and there exists a  sequence  $(\vphi_{j})$ of \psh functions
in $\mcal E_{0} (\W)$ such that $\vphi_{j} \searrow \vphi$ 
and $ \lim_{j \to + \infty} \int_{\W} (-\vphi_{j})^p 
(dd^c \vphi_{j})^n  =  \int_{\W} (-\vphi)^p (dd^c \vphi)^n.$ 
\end{lem} 
 Such a sequence will be called a $p-$admissible sequence decreasing to $\vphi.$
\bigskip
\proof Observe that for any  sequence 
$(\vphi_{j})$ of \psh functions in $\mcal E_{0},$ decreasing
 to $\vphi,$ with the condition (\ref{eq:Ep}), the sequence
$(- \vphi_{j})^p$ is an increasing sequence of 
lower semi-continuous  functions on $\W$ converging to 
$(- \vphi)^p$ and the  sequence of measures $(dd^c \vphi_{j})^n$ 
 converges weakly to $(dd^c \vphi)^n$ on $\W.$ Therefore it 
  follows that $ \int_{\W} (- \vphi)^p (dd^c \vphi)^n \leq 
  \liminf_{j \to + \infty} \int_{\W} (- \vphi_{j})^p 
  (dd^c \vphi_{j})^n,$ which proves that 
  $ e_{p} (\vphi) < + \infty.$
For $p = 0,$ the result of the lemma  follows from the definition.
 Now fix  $p > 0.$ Then, since for the given function 
$\vphi \in \mcal E_{p}$ the measure $(dd^c \vphi)^n$ 
puts no mass on pluripolar sets, it follows from  Cegrell's 
decomposition theorem [Ce 2] that there exists 
 $ \psi_{0} \in \mcal E_{0} (\W)$ and $ 0 \leq f \in L^{1} 
(\W;(dd^c \psi_{0})^n)$ such that 
$(dd^c \vphi)^n = f (dd^c \psi_{0})^n.$ By 
Kolodziej's theorem (see [Ko 1], [Ce 2]),  for any integer 
$j \geq 1$
there exists  $\vphi_{j} \in \mcal E_{0} (\W)$ such that 
$(dd^c \vphi_{j})^n = \min \{f , j \} \cdot (dd^c \psi_{0})^n.$ Now by the 
comparison principle [Be-Ta 1], we see that the sequence 
$(\vphi_{j})$ decreases to a \psh function $\psi$ 
on $\W$ such that $\psi \in \mcal E_{p} (\W)$ 
and $(dd^c \psi)^n = (dd^c \vphi)^n.$ 
By the comparison principle
it follows that $\vphi = \psi$ on $\W.$ Now by the monotone 
convergence theorem we obtain 
 $ \int_{\W} (- \vphi)^p (dd^c \vphi)^n = \lim_{j \to + \infty} 
 \int_{\W} (- \vphi_{j})^p  \min (f,j) (dd^c \psi_{0})^n = 
 \lim_{j \to + \infty} \int_{\W} (- \vphi_{j})^p (dd^c \vphi_{j})^n,$
which proves the lemma.
 \fin 

  Now let us give a quantitative characterization of the class 
  $\mcal E (\W)$ in terms of some capacity introduced by Bedford (see [Be]).
  
Given  $\vphi \in PSH^{-} (\W)$ and $K \sub \W$ a borelean 
set, following [Be] we define the following positive "$\vphi-$capacity"  
  $$ C_{\vphi} (K; \W) := \sup \{ \int_{K} (dd^c \psi)^n \ ; \ \psi \in 
 PSH^{-} (\W) \cap L^{\infty} (\W), \vphi \leq 
 \psi \leq 0 \}.$$
  
  Consider also the corresponding extremal function
  \begin{equation}
        \Tilde{\vphi}_{K}  = \sup \{u \in PSH^- (\W) ;  u \leq \vphi \ 
 \mathrm{  q.e. \  on \  } \ \ K \}.
        \label{eq:tildephi}
  \end{equation}
  where $u \leq \vphi$ q.e. (quasi-everywhere)  on $K$ means
  outside a pluripolar subset of $K.$
 
 Then we get the following characterization of functions from the class 
 $\mcal E (\W).$
\begin{prop} Let $\vphi \in PSH^- (\W).$  Then 
$\vphi \in \mcal E (\W)$ if and only if for any 
compact set $K \sub \W$ we have  $C_{\vphi} (K) < + \infty.$ 
 Moreover, if $\vphi \in \mcal E (\W),$ then for any
 borelean set $K \Subset \W,\Tilde{\vphi}_{K} \in 
\mcal F (\W),\Tilde{\vphi}_{K} = \vphi$ q.e. on  $K$ 
and 
\begin{equation}
        C_{\vphi} ({K}^{\circ}) \leq \int_{\W} (dd^c 
 \Tilde{\vphi}_{K})^n \leq C_{\vphi} (\overline{K};\W).
        \label{eq:Cineq}
\end{equation}
 \end{prop} 
  \proof  From the 
  definition of $\Tilde{\vphi}_{K}, $ it follows that
  $\Tilde{\vphi}_{K}$ is \psh on $\W,$ 
  maximal on $\W \setminus \overline {K}$ and satisfies
the inequality  $ \vphi \leq \Tilde{\vphi}_{K} \leq 0$ on 
$\W.$ \\
 Now assume that  $\vphi \in \mcal E (\W).$ Then 
  it follows  from ([Ce 2], [Ce 3]) that $\Tilde{\vphi}_{K}  = \sup\{\vphi, \Tilde{\vphi}_{K} 
  \} \in \mcal E (\W)$ and $ \int_{\W} 
  (dd^c \Tilde{\vphi}_{K} )^n = \int_{\overline{K}} 
  (dd^c \Tilde{\vphi}_{K} )^n  < + \infty.$ 
  
 Let us first prove  that 
$\Tilde{\vphi}_{K} \in \mcal F (\W).$   Indeed, take a decreasing 
sequence $(\vphi_{j})_{jÊ\geq 0}$ from $\mcal E_{0} (\W)$ 
converging to $\vphi$ on a neighbourhood of $\overline{K}$ such that 
$\sup_{j} \int_{\W} (dd^c \vphi_{j})^n < + \infty.$  Then 
for each $j \in \N$ define the function $\Tilde {\vphi}_{j}$ by 
the formula (\ref{eq:tildephi}) with $\vphi$ repalced by 
$\vphi_{j}.$ Then $(\Tilde {\vphi}_{j})$ is a decreasing 
sequence of \psh  functions from the class $\mcal E_{0} (\W)$  
such that  ${\vphi}_{j} \leq \Tilde {\vphi}_{j}$ on $\W$ and 
$\Tilde {\vphi}_{j} = \vphi_{j}$ q.e. on $K.$ Therefore 
$(\Tilde {\vphi}_{j})$ converges to a \psh function $\psi$ such that
$\Tilde{\vphi}_{K} \leq \psi $ on $\W$ and $\psi = \vphi$ q.e. on 
$K.$ Thus $\psi =  \Tilde{\vphi}_{K}.$ Since ${\vphi}_{j} \leq 
\Tilde {\vphi}_{j}$ on $\W$ and these functions belong
 to $\mcal E_{0} (\W),$ it 
 follows that $\int_{\W} (dd^c \Tilde {\vphi}_{j})^n \leq 
 \int_{\W} (dd^c \vphi_{j})^n$ for any $j\geq 0.$ Thus the sequence 
 $(\Tilde {\vphi}_{j})_{j \geq 0}$ decreases to 
 $\Tilde {\vphi}_{K}$ and  $\sup_{j} \int_{\W} (dd^c \Tilde 
 {\vphi}_{j})^n \leq \sup_{j} \int_{\W} (dd^c \vphi_{j})^n < + 
 \infty,$ which proves that  $\Tilde{\vphi}_{K} \in \mcal F (\W).$  
 Then since $\vphi \leq \Tilde 
 {\vphi}_{j},$ it follows that
  $\int_{\W} (dd^c \Tilde {\vphi}_{j})^n = \int_{K} 
  (dd^c \Tilde {\vphi}_{j})^n  \leq C_{\vphi} (K).$ 
 By Cegrell (see [Ce 3]), the sequence of measures
 $(dd^c \Tilde {\vphi}_{j})$ 
 converges to the measure $(dd^c \Tilde {\vphi})^n,$  
 thus $\int_{\W} (dd^c \Tilde {\vphi}_{K})^n  \leq  \liminf_{j} 
 \int_{\W} (dd^c \Tilde {\vphi}_{j})^n  \leq C_{\vphi} 
 (\overline{K})$, which 
 proves the second inequality in (\ref{eq:Cineq}). \\ 
  Now let $\psi \in \mcal E_{0} (\W)$ be choosen so that 
$ \vphi \leq \psi$ on $\W$ and set $\psi_{j} := \sup 
\{\psi,\vphi_{j}\}$ on $\W.$ Then $\vphi_{j} \leq \psi_{j}$
 and then $\Tilde{\vphi}_{j} \leq \psi_{j}$ q.e. on $K.$
 Since functions from $\mcal E_{0} (\W)$ put no mass on 
 pluripolar sets, it follows from Demailly's inequality
 ([De]) that 
  \begin{eqnarray}
  \int_{K} (dd^c \psi_{j})^n  & \leq & \int_{K} 
        (dd^c \sup \{\psi_{j} , \Tilde{\vphi}_{j}\})^n
          \\
         & \leq & \int_{\W} 
        (dd^c \sup \{\psi_{j} , \Tilde{\vphi_{j}}\})^n
        \nonumber  \\
         & \leq & \int_{\W} (dd^c \Tilde{\vphi_{j}})^n = 
         \int_{\overline{K}} (dd^c \Tilde{\vphi_{j}})^n, 
         \nonumber \\
 \end{eqnarray}                   
 where the last inequality follows from the comparison principle for 
 functions in $\mcal E_{0} (\W).$
 
  Therefore by the convergence theorem (see [Ce 2], [Ce 3]) we have
  $$ C_{\vphi} ({K}^{\circ};\W) \leq \int_{\W} 
         (dd^c \Tilde{\vphi}_{ K})^n < + \infty,$$
 which proves the first inequality in (\ref{eq:Cineq}) as well as the 
 necessary condition of the theorem.\\
 Now assume that the condition on the  capacity $C_{\vphi}$ is 
satisfied. Then consider a decreasing sequence of \psh  functions 
$(\vphi_{j})$ from the class $\mcal E_{0} (\W)$ converging to 
$\vphi.$ Take any subset $K \Subset \W$  and define $
\Tilde{\vphi}_{j} := (\Tilde{\vphi}_{j})_{K}$ for $j \in \N^{*}.$
Then $\Tilde{\vphi}_{j} \in \mcal E_{0} (\W)$ and 
$\sup_{j} \int_{\W} (dd^c \Tilde{\vphi}_{j})^n = \sup_{j} 
\int_{\overline{K}} (dd^c \Tilde{\vphi}_{j})^n 
\leq C_{\vphi} (\overline{K};\W) < + \infty.$
We know from the first part that  $(\Tilde{\vphi}_{j})$ decreases 
to $\Tilde{\vphi}_{K} = \vphi$ q.e. on $K.$ Therefore we have 
proved  that  $\vphi \in \mcal E  (\W).$ Moreover by the convergence theorem
 we obtain the second inequality in (\ref{eq:Cineq}). 
 \fin 
 Bedford considered the following class (see [Be]).   Let 
  $\theta  : \R \longrightarrow \R$ be a monotone decreasing 
  function such that 
  \begin{equation}
        \int_{1}^{+ \infty} \frac{\theta  (t)}{t} d t
         < + \infty
         \label{eq:growth}
  \end{equation}
  and the function $t \longmapsto - (- t \ \theta (- t))^{1 \slash n}$
  is monotone increasing and convex on $]- \infty , 0[.$
  Then define  $\mcal B (\W) $ to be the class of negative 
  function $\psi \in PSH (\W) $ such that for any $z_{0} 
  \in \W$ there  exists a neighbourhood $\omega$ of 
  $z_{0},$  a negative 
  \psh function $v$ on $\omega$ and  a decreasing function 
  $\theta$ satisfying (\ref{eq:growth}) such that   
  $ - (- v \theta (- v))^{1 \slash n} \leq \psi $ on 
  $\omega.$ \\
 From the last result we can deduce the following one which 
 provides concrete examples of functions from the class $\mcal E (\W).$ This
  result has been also obtained by the first author (see [Ce 4]).
 \begin{prop} For any hyperconvex domain $\W \Sub \C^n,$ we have
 $\mcal B (\W) \sub  \mcal E (\W).$ In particular, for any negative \psh 
 function $u $ on $\W$ and any 
 $0 < \alpha < 1\slash n,$ $ - (- u)^{\alpha} \in \mcal E (\W).$
 \end{prop}
 \proof Bedford has proved that for any function $\psi \in \mcal 
 B (\W), $ the condition $C_{\psi} (K;\W) < + 
 \infty$ holds for any $K \Subset \W$ (see [Be]). 
 Therefore the inclusion 
 $\mcal B (\W) \sub  \mcal E (\W)$ follows from
 the last proposition. Since $\alpha < 1\slash n,$ we clearly have 
 $- (- u)^{\alpha} =  - (- u \theta (- u))^{1 \slash n},$ where $ 
 \theta (t) = t^{n \alpha - 1}$ which satisfies the condition (\ref 
 {eq:growth}). Therefore $ - (- u)^{\alpha} \in \mcal B 
 (\W) \sub  \mcal E (\W).$
 \fin

 \section{Capacity of sublevel sets of \psh  \\
 functions in subclasses of $\mcal E (\W)$}    
 
 Now we prove the following capacity estimate of the sublevel sets of 
 \psh functions of finite energy.
  \begin{prop} Let $\vphi \in \mcal E_{p} (\W)$ if $p > 0$ and 
  $\vphi \in \mcal F (\W)$ if $p = 0.$  Then the following estimate
  \begin{equation}
        \mathrm{Cap} (\{ z \in \W \ ; \ \vphi (z)  < - s \};\W) \leq 
  c_{n,p} \cdot  e_{p} (\vphi) \cdot s^{- n - p} , \forall s > 0,
        \label{eq:capest}
 \end{equation}
 holds, where $c_{n,p} > 0$ is an absolute constant and 
  $e_{p} (\vphi) $ is the $p-$energy of $\vphi.$ 
 \end{prop}
  \proof 1) Assume first that $\vphi \in \mcal E_{0} (\W)$ 
 and define $\W (\vphi,s) := \{z \in \W ; \vphi (z) < - s \}$ for 
 each $s > 0.$  
 Let $K \sub \W (\vphi,s) $ a fixed pluriregular subset and 
 $h_{K}$ the $(-1,0)-$extremal function of the condenser  
 $(K,\W).$  Then $h_{K} \in \mcal E_{0} (\W)$ and since 
 $- \vphi \slash s \geq 1$ on $K,$ we have
 \begin{eqnarray}
        \mathrm{Cap} (K;\W) & = & \int_{K} (dd^c h_{K})^n \leq \int_{K} 
 \bigl(\frac {- \vphi}{s}\bigr)^{n + p} (dd^c h_{K})^n
        \label{}  \\
                   &    \leq &  \frac{1}{s^{n + p}} \int_{\W} 
                   (-\vphi)^{n + p}  (dd^c h_{K})^n
        \nonumber  \\
         & = & \frac{c_{n,p}}{s^{n + p}} \int_{\W} (-\vphi)^p 
         (dd^c\vphi)^n,
        \nonumber
 \end{eqnarray}
 where the last inequality follows by integration by parts 
 (see [Ce 3], [Bl]).
 Then we deduce that
 $$ \mathrm{Cap} \bigl(\W (\vphi, s) ; \W\big) \leq 
 \frac{c_{n,p}}{s^{n + p}} 
 \int_{\W} (-\vphi)^p  (dd^c\vphi)^n$$
  for each $\vphi \in \mcal E_{0} (\W).$ \\
2)  Let $\vphi \in \mcal E_{p} (\W)$ and $(\vphi_{j})$  
a $p-$admissible sequence of \psh functions  from the class 
$\mcal E_{0} (\W)$ decreasing to $\vphi.$  Then applying the  
estimate (\ref{eq:capest}) to each function $\vphi_{j} \in \mcal E_{0} (\W)$ 
we obtain the following estimate
  $$ \mathrm{Cap} \bigl(\W (\vphi_{j},s) ;\W\bigr) \leq 
  \frac{c_{n,p}}{s^{n + p}} 
  \int_{\W} (-\vphi_{j})^p (dd^c \vphi_{j})^n, \forall j .$$
 Now by Lemma 2.1, there exists a decreasing sequence of \psh 
functions  $(\vphi_{j})$  from the class $\mcal E_{0} (\W)$ 
which converges to $\vphi$ and satisfies the condition 
$e_{p} (\vphi) = \lim_{j} \int_{\W} (- \vphi_{j})^p 
(dd^c \vphi_{j})^n.$ Then applying the last estimate to this 
function we get the following estimate
  $$\mathrm{Cap} \bigl(\W (\vphi,s);\W) \leq  c_{n,p} \frac{e_{p}  
 (\vphi)}{s^{n + p}}.$$
This proves our proposition.   \fin 
 Let us prove a converse to the last result which shows that the estimates
 (\ref{eq:capest}) are almost sharp.
\begin{prop} Let $\vphi \in \mcal E (\W)$ a function such that there 
 exists an open subset $\w \subset \W,$ a constant $A > 0$ and a 
 real number $q > n$ such that 
 $$ \mathrm{Cap} (\{z \in \w \ ; \  \vphi (z) < -s \}; \W)  \leq 
 C s^{- q}, \forall s > 0.$$
Then $\Tilde{\vphi}_{\w} \in \mcal E_{p} (\W)$ for any real 
number $p$ with $0 < p < q - n.$ 
 \end{prop}
 \proof First we claim that for any \psh function
  $u \in \E (\W)$ and 
any Borel set $B \subset \W,$
we have 
 \begin{equation}
\int_{B} (dd^c u)^n \leq  (\sup_{B} 
\vert u \vert)^n \ \mathrm{Cap} (B;\W),
        \label{eq:m-est}
 \end{equation}
provided that $\sup_{B} 
\vert u \vert < + \infty.$ \\
To prove this estimate, set $M := \sup_{B} 
\vert u \vert $ and 
define the function $v = \sup \{ u \slash M  , - 1\}$ on $\W.$ 
Then $v \in \mcal PSH (\W), - 1 \leq v \leq 0$ and 
$ v = u \slash M $ on $B.$ 
 Therefore from Demailly's inequality ([De]),  it follows  that 
 $ (M ^{- n} \int_{B}(dd^c u)^n  \leq 
 \int_{B} (dd^c v)^n 
 \leq \mathrm{Cap} (B;\W),$ 
 which proves our claim.
 
 Now define the following sets $B_{j} := \{z \in  {\w} \ ; \ - 2^{j + 
 1}  \leq \vphi (z)  < - 2^j \}$ for $j \geq 0.$
 Then it follows from (\ref{eq:m-est})  that
 $$\int_{B_{j}} (- \vphi)^p (dd^c \vphi)^n \leq  (\sup_{A_{j}} 
 \vert \vphi \vert)^{n + p} \mathrm{Cap} (B_{j}; \W) \leq  C 2^{j (n + p - 
 q)} $$
 and then $\int_{\w} (- \vphi)^p (dd^c \vphi)^n \leq C
  \sum_{j \geq 0} 
 2^{(n + p - q) j} < + \infty,$ since $p < q - n.$ Therefore we have
 $\int_{\W} (- \Tilde{\vphi}_{\w})^p (dd^c \Tilde{\vphi}_
 {\w})^n < + \infty.$
 \fin
 As a consequence we state the following result which improves the result of
 Proposition 2.3 and provides examples of
 functions in the classes $\mcal E_p (\W).$
 \begin{cor} Let $u \in PSH^{-} (\W),$ $\alpha $ a real number such 
 that $0 < \alpha  < 1 \slash n$ and $\vphi = \vphi_{\alpha} :=  
 - (- u)^{\alpha}.$ Then for any Borel  subset $\w \Subset 
 \W$ $  \Tilde{\vphi }_{\w} \in \mcal 
  E_{p} (\W)$ for any real number $p$ such 
 that $0 < p < 1 \slash \alpha - n,$  
 \end{cor} 
\proof  Indeed it is easy to check that $\mathrm{Cap} (\{z \in \w \ ; \ 
 \vphi (z) < - 
 s\} ; \W) \leq A s^{- 1 \slash \alpha},$ for any 
 $s > 0.$ Thus the result follows from the last 
 one. \fin 
 
It is possible to characterize the class $\mcal F^a (\W)$ by means of the
behaviour of the  capacity
of sublevel sets.
\begin{prop} Let $\vphi \in \mcal F (\W).$ Then the following properties are
  equivalent\\
$(i)$ $\vphi \in \mcal F^a (\W),$ \\
$(ii)$ $\int_{\{\vphi = - \infty\}} (dd^c \vphi)^n = 0,$ \\
 $(iii)$ $ \lim_{s \to + \infty} s^{n} \mathrm{cap} (\{ \vphi < - s\}; \W) = 0.$
\end{prop}
\proof Take a sequence  $(\vphi_j)$ of continuous functions from  $\E_0 (\W)$ which decreases
to  $\vphi$ and satisfies $\sup_j \int_\W (dd^c \vphi_j)^n < + \infty.$

Define the open sets  $\W_j (s) := \{\vphi_j < - s \}$, $\W (s)
:=  \{\vphi < - s \}$ and the functions
$$ a_j (s) := \mathrm{cap} (\{ \vphi_j < - s\}; \W), \ \ a (s) := \mathrm{cap} (\{
\vphi < - s\}; \W),$$
and 
$$ b_j (s) := \int_{\W_j (s)} (dd^c \vphi_j)^n , \ \ \ b (s) :=  \int_{\W (s)}
(dd^c \vphi)^n.$$
We claim that for $j \in \N$ and $s > 0,$
\begin{equation}
s^n a_j (2 s) \leq b_j (s)  \leq s^n a_j (s).
\label{eq:compest}
 \end{equation}
Indeed take any function $u \in PSH (\W)$ with
 $- 1 \leq u \leq 0.$ Then $\W_j (2 s) \sub \{\vphi_j
 \slash s < u - 1\} \sub \W_j (s) \Sub \W.$ By the comparaison
principle, we get
$$\int_{\W_j (2 s)} (dd^c u)^n \leq 
\int_{ \{ \vphi_j \slash s < u - 1\}} (dd^c  \vphi_j)^n \leq 
\int_{\W_j (s) } s^{- n } (dd^c  \vphi_j)^n = s^{- n} b_j (s).$$
Taking the supremum over all $u$'s, we obtain the first inequality of
(\ref{eq:compest}).

  To obtain the second inequality, observe that for $0 < s < t,$ 
$\sup\{\vphi_j , - t\} = \vphi_j $ on the open set
 $ \{\vphi_j >  - t \}$
which is a neighbourhood of $\partial \W_j (s)$ and then
$\int_{ \W_j (s)} (dd^c \sup\{\vphi_j , - t\})^n = 
\int_{ \W_j (s)} (dd^c \vphi_j)^n.$
Therefore 
$$a_j (s) \geq t^{- n} \int_{ \W_j (s)} (dd^c \sup\{ \vphi_j , - t\})^n =  t^{-
  n} \int_{ \W_j (s)} (dd^c  \vphi_j)^n,$$
which  proves the required inequality since $t>s$ is
 arbitrarily close to $s.$
Taking the limit in  (\ref{eq:compest}) when $j \to + \infty,$ we obtain
\begin{equation}
s^n a (2 s) \leq b (s^-) \ \ \mathrm{  and  } \ \ b (s) \leq s^n a (s), \
\forall s > 0.
\label{eq:compest2}
 \end{equation}
where $b (s^-) := \lim_{\vep \to 0^+} b (s - \vep).$

From the estimates (\ref{eq:compest2}), it follows that the conditions $(ii)$
 and $(iii)$ are equivalent. Moreover it is clear that $(i)$ implies
 $(ii)$. So it is enough to prove that $(iii) $ implies $(i).$
Indeed, assume that that $\lim_{s \to + \infty} s^{n} a
(s) = 0$ and take a pluripolar subset $K$ of $\W$.
It follows from ([De], [Ce-Ko]) that
$$ \int_{K \setminus \W (s)} (dd^c \vphi)^n \leq \int_{K \setminus \W (s)}  (dd^c
\sup\{\vphi , - s\})^n  \leq s^{n} \mathrm{cap} (K;\W) = 0. $$
Moreover, by  (\ref{eq:compest2}), we have
$$\int_{K \cap \W (s)} (dd^c \vphi)^n \leq b (s) \leq s^{n} a (s).$$
Therefore $\int_K (dd^c \vphi)^n = 0$ and then $\vphi \in \F^a (\W).$ \fin 
 \section {Global subextension of \psh \\
functions with weak singularities} 
 Here we want to prove a general subextension theorem for a class 
 of \psh functions of weak singularities, generalizing a theorem by 
 El Mir ([El])  and also by Alexander and Taylor (see [Al-Ta], [De]).
 Then we will apply our result to derive
 theorems on subextension of \psh functions of 
 finite energy.

To sate our results we need to introduce the usual Lelong classes of
\psh functions.
\begin{equation}
\mcal L_\gamma (\C^n) := \{u \in PSH (\C^n) ; \limsup_{r \to + \infty}\frac{\max_{\vert
  z\vert = r} u (z)}{\log r} \leq \gamma\}, \gamma > 0.
\end{equation}
When $\gamma = 1$ we write $\mcal L (\C^n) = \mcal L_1 (\C^n).$

 \begin{thm} Let $\vphi \in  PSH^{-} (\W)$ and $\w \sub \W$ an open 
 subset. Define the function $ \chi (s) = \chi_{\vphi} (s,\w) := \mathrm{Cap} (\{z \in 
 \w \ ; \ \vphi (z) < - s \}; \W).$  Assume that the following 
 integral condition
 \begin{equation}
 \int_{1}^{+ \infty} \chi (s)^{1 \slash n} d s < + 
 \infty 
        \label{eq:cap-cond}
 \end{equation}
 holds. 
  Then for any $\vep > 0,$ there exists a function $U_{\vep} \in \mcal
  L_{\vep} (\C^n)$ such that $U_{\vep} \leq \vphi$ on  $\w.$
 
 In particular  $\nu_{\vphi} (a) = 0$ for any $a \in \w.$
 \end{thm}
\proof We use the same construction as in  ([Al-Ta]).
 Indeed let us denote by 
 $$ M (s) := \max_{\overline \W} V_{s} = - \log  
 T_{\overline \W} \bigl(\w (\vphi;s)\bigr), \ s > 0,$$
 where   $V_{s}$ is the $\mcal L-$extremal function of the open 
 set  $\w (\vphi;s) := \{z \in \w \ ; \ \vphi (z) < - s\}.$ Then by
Alexander-Taylor's inequality ([Al-Ta]), we deduce the 
following
  \begin{equation}
   M (s) \geq  \chi (s)^{- 1 \slash n}, \ \
 \forall s > 0,
        \label{eq:est-M}
  \end{equation}
  where $\chi (s) = \chi_{\vphi} (s;\w)$ for $s > 0.$
 Now define the following function
 \begin{equation}
    w_{s} (z) = w (z,s)  := V_{s} (z) - M (s) , \  
        (z,s) \in \C^n \x \R^{+}.
        \label{eq:ws}
 \end{equation}
 We claim that this function satisfies the following  properties
 \newline $(i)$ $w_{s} \in \mcal L$ and $w_{s}  (z) \leq a +  
 \log^{+} \vert z \vert , \ 
 \forall z \in \C^n, \ \forall s > 0,$
 \newline $(ii)$ $w_{s} (z) = - M (s), \ \forall z \in \w (\vphi,s), 
 \forall s > 0,$
 \newline $(iii)$ $\max_{\overline \W} w_{s} = 0$ and 
 $\int_{\W} w_{s} (z) d \lambda (z) \geq - b, \  \forall s > 0,$ \\
 where $a, b > 0$ are absolute constants. 
 
 Assume for the moment that all the above properties are 
 satisfied and observe that for any fixed $z \in \C^n,$ the function
 $s \longmapsto w_{s} (z)$ is a function of bounded variation( equal to the difference of two 
 monotone functions) and upper bounded on $\R^{+},$  by 
 condition $(i).$ Therefore we can define the following function
 $$  v_{c} (z) :=  \int_{c}^{+ \infty} w (z,s) \chi (s)^{1 \slash n} 
 d s, \ z \in \C^n,$$
 for each $c > 0.$  From condition $(i)$, it follows that
 \begin{equation}
        v_{c} (z) \leq a \cdot  \eta_{c} + \eta_{c} \cdot \log ^{+} 
        \vert z \vert, \ \forall z \in \C^n,
        \label{eq:est-upper}
 \end{equation}
 where $\eta_{c} := \int_{c}^{+ \infty} \chi (s)^{1 \slash n} d s.$
 Now from  $(iii)$ it follows that
 \begin{equation}
        \int_{\W} v_{c} (z) d \lambda_{2 n} (z) \geq - b \cdot \
        \eta _{c}
        \label{eq:int}
 \end{equation}
 Then  it follows from (\ref{eq:int})
 and (\ref{eq:est-upper}) that $ v_{c}$ is \psh on $\C^n.$
 Now fix $ t > c \geq 0$ and $z \in \w (\vphi,t).$ Then  by $(ii)$, 
 for any $ s < t,$  we have $ \vphi (z) < - t$ and $w_{s} (z) = - 
 M (s).$ Therefore, 
 since $w_{s} \leq 0 $ on $\W,$  we get from (\ref {eq:est-M}) 
 the following estimate
  $$ v_{c} (z) \leq \int_{c}^t w_{s} (z) \chi (s)^{1 \slash n} 
  d s \leq ( - t + c).$$ 
  This means that $v_{c} (z) \leq  
  \vphi (z) + c$ if $ \vphi (z) < - c.$ But if $\vphi (z) \geq - c,$ 
  this inequality is clearly satisfied, since $ v_{c} (z) \leq 0$ 
  for any $z \in \W.$
  Define  $u_{c} (z) :=  v_{c} (z)  - c, $ for $ z \in \C^n.$ 
  Then it is clear that $  u_{c} \leq \vphi$ on 
  $\w.$ Moreover, from (\ref{eq:est-upper}), it follows that  
  \begin{equation}
         u_{c} (z) \leq a  \cdot \eta _{c} - c +  \eta_{c}  
        \cdot \log ^{+} \vert z \vert, \ \forall z \in \C^n,
        \label{eq:upper-est}
\end{equation}

Now given $\vep > 0,$ we can choose $c = c (\vep) > 0$ such that
$\eta (c) < \vep,$ then corresponding function $U_{\vep} := u_{c (\vep)}$ 
satisfies the conclusions of the theorem with $\gamma (\vep) :=  a \cdot \eta
_{c (\vep)} - c (\vep).$\\

Now it remains to prove that our function (\ref {eq:ws}) satifies 
the properties $(i), (ii)$ and $(iii).$
By definition $w_{s} \in \mcal L$ and $\max_{\overline \W} 
w_{s} = 0.$  Then $w_{s} \leq V_{\overline \W}$ on $\C^n$ for 
any $s > 0,$ which proves $(i),$ since $  V_{\overline \W} \in 
\mcal L.$ 
The condition $(ii)$ is trivial since 
$V_{s} = V_{\w (\vphi,s)} = 0$ on $\w (\vphi,s).$ The  
condition $(iii)$ is related to an inequality by Alexander 
(see [Al], [Sic 2], [De]) and can be proved easily as follows. 
Observe that the normalized 
subclass $\Dot{\mcal L}_{\overline \W} := \{ w \in \mcal L \ ; \ 
\max_{\overline \W} = 0 \}$ is a (relatively) compact subset of 
$\mcal L$ for the $L_{loc}^{1}-$topology (see [Ze]) and the functionnal $w \longmapsto 
\int_{\W} w  (z) d \lambda (z) $ is continuous on $\mcal L.$ 
Therefore it is bounded on $\Dot{\mcal L}_{\overline \W},$ which 
proves the condition $(iii).$
\fin 

 Now from our result we can deduce the  Alexander-Taylor's
 subextension  theorem. 

 Let $h : \R^- \longrightarrow \R^-$ be an 
increasing convex function such that
  \begin{equation}
  \int_{1}^{+ \infty} \frac{- h (-t)}{t^{1 + 1\slash n}} dt < + 
  \infty.
        \label{eq:h}
  \end{equation}
  Then we obtain the following result.
  \begin{cor} Let $u \in PSH^{-} (\W)$ and $h : \R^- \longrightarrow \R^-$ be
 an increasing convex function satisfying the condition (\ref{eq:h}). 
Then for any subdomain $\w \Sub \W,$ for any $\vep > 0,$ there 
exists a function $U_{\vep} \in \mcal L_{\vep} (\C^n)$ such that $U_{\vep} \leq  h (u) $ on $\w.$

 \end{cor}
 \proof Let $g : \R^- \longrightarrow \R^-$ be the inverse function 
 of $h.$ Then $ \w (h (u);s) = \w (u;g(- s))$ for any $s > 0.$
 Now use the usual capacity estimate for $u$ to conclude that
 \begin{equation}
     \mathrm{Cap} (\w (h (u);s); \W) \leq \frac{A}{- g(- s) } , \forall s > 0.
        \label{eq:caph}
 \end{equation}
 Now observe that the condition  (\ref{eq:h} ) on $h$ implies that
 $ \int_{1}^{+ \infty} (- g(- s))^{- 1 \slash n} ds <  + \infty $. Therefore
 from the estimate (\ref{eq:caph}), it follows that
 the condition  (\ref{eq:cap-cond})  is 
 satisfied for the function $h (u)$ and then the corollary follows from the
 last theorem.
 \fin 

 Now using capacity estimates from section 3 and the last theorem, we easily see that functions from the 
 classes  $\mcal E_{p} (\W),$ with  $p > 0,$  have global subextension of
 arbitrary small logarithmic growth at infinity.
\begin{cor} Let $\vphi \in \mcal E_{p} (\W),$ with $p > 0.$ 
Then for any $\vep > 0,$ there exists a
function  $U_{\vep} \in \mcal L_{\vep} (\C^n)$ such that
 $U_{\vep} \leq \vphi$ on $\W.$
 \end{cor}
 \proof From the estimates  (\ref{eq:capest})  of Proposition 3.1,
it follows that the condition (\ref{eq:cap-cond}) of Theorem 4.1 
is satisfied with $\w = \W,$ which implies our result.
\fin 

 \section {Global subextension of psh functions with \\ uniformly bounded
   Monge-Amp\`ere masses} 
  As we pointed out in the introduction, on any smoothly bounded 
 domain in $\C^{2}$ there is a smooth \psh function which admits no
 subextention to any larger domain 
 (see [Be-Ta 3]).
 In contrast to this negative result,the first and the third authors  proved
 that 
 for any hyperconvex domain $\W \Sub \C^n,$ functions from
 the 
 class $\mcal F (\W)$ always admit a subextension to any
 larger bounded hyperconvex domain (see [Ce-Ze]).

Here we want to prove that such functions have a global subextension which is \psh
 of logarithmic growth on $\C^n.$

Besides the Lelong classes $\mcal L_{\gamma} (\C^n)$ defined in
section 4, we also need the following class.
$$\mcal L_{\gamma}^+ (\C^n) := \{u \in PSH (\C^n) ;
\sup_{z \in \C^n} \vert u (z) - \log^{+} \vert z \vert \vert < + \infty\}.$$

Now we can state our main result.
 \begin{thm} Let $\W \Subset \C^n$ be a bounded
 hyperconvex domain and $\vphi \in \mcal F (\W).$
 Then  there exists a \psh function $u \in 
 \mcal L_{\gamma} (\C^n),$ with $\gamma^n := \int_\W
 (dd^c \vphi)^n$ such that
 $\max_{\ove{\W}} u = 0$ and $u \leq \vphi$
 on $\W$. 
 \end{thm}
 \proof
1) Assume first that $\vphi \in \mcal E_0 (\W)$ and define
 the following Borel measure $\mu := {\bf 1}_\W (dd^c \vphi)^n.$
 Fix a ball $\B \sub \C^n$ such that $\ove \W \sub \B.$
 Then in general there is no bounded \psh function $v$
 on $\B$ such that $ (dd^c v)^n \geq \mu$ on $\B$ (see the
 example below).
  We will approximate the measure $\mu$ by measures for which
 such bounded \psh functions exist.
 Indeed, since $\mu$ puts no mass on pluripolar sets, by 
 ([Ce 1]) there exists $\psi \in \mcal E_0 (\B)$ and
 $f \in L^1 (\B , \mu)$ such that $\mu = f \cdot
 (dd^c \psi)^n$ on $\B.$ Then consider the sequence of
 measures $\mu_k := {\bf 1}_{\W} \inf \{f,k\} (dd^c \psi)^n,
 k \in \N$ with compact support in $\B.$
 
 Fix an integer $k \geq 1.$ Since $\mu_k \leq (dd^c \psi_k)^n$ on $\B,$ where $\psi_k
:= k^{1 \slash n} \psi \in \mcal E_0 (\B),$ it follows from
([Ko 2]) that there exists $u_k \in \mcal L_{\gamma_k}^+ (\C^n)$
 such that $(dd^c u_k)^n = \mu_k$ on $\C^n,$ where $\gamma_k^n := \mu_k (\B).$
We can normalize $u_k$ so that
$\max_{\ove{\W}} u_k = 0.$
We can also find $g_k \in \mcal E_0 (\W)$ such that
$(dd^c g_k)^n = \mu_k$ on $\W.$
Then from the comparison principle, we have
$u_k \leq g_k$ on $\W$ and since the sequence of measures $(\mu_k)$ is
increasing, the sequence of \psh functions $(g_k)$ decreases to $\vphi$ on $\W$.
By Hartogs lemma, $u := (\limsup_{k \to + \infty} u_k)^*$ is plurisubharmonic
on $\C^n$ and $\max_{\ove{\W}} u = 0.$ It is clear that
$u \leq \vphi$ on $\W$ and $u \in \mcal L_{\gamma} (\C^n),$ where $\gamma^n
:=\int_{\ove \W} (dd^c u)^n,$ since  $\gamma_k \leq \gamma,
\forall k \in \N.$

 2) Assume now that $\vphi \in \mcal F (\W).$ By Lemma 2.1, there exists a decreasing sequence 
 $(\vphi_{j})$ of functions from the class $\mcal E_{0} (\W)$
 which converges to $\vphi$ on $\W$ and $\int_{\W}
 (dd^c \vphi)^n = \lim_{j} \int_{\W} (dd^c \vphi_{j})^n.$
  Let us define $ \gamma > 0$ so that
 $\gamma^n := \int_{\W} (dd^c \vphi)^n < + \infty$ and
fix  $j \in \N.$ Then by the first case there exists $u_j \in \mcal L_{\gamma_j}^+
 (\C^n)$ such that $\max_{\ove \W} u_j = 0$ and
$u_j \leq \vphi_j$ on $\W,$ where $\gamma_j^n = \int_\W (dd^c \vphi_j)^n.$

Again the function $u := (\limsup_{j \to + \infty} u_j)^* \in \mcal L_{\gamma}
(\C^n)$ and satisfies the inequality $u \leq \vphi$ on $\W$
and by Hartogs' lemma we have $\max_{\ove \W} u = 0.$
\fin 

 From this result we get the following one.
 \begin{cor} Let  $\W  \Subset \C^n$ 
 be a bounded hyperconvex domain and $\vphi \in \mcal E (\W).$ Then 
 for any open set  $\w \Subset \W,$ there exists a function 
 $u \in \mcal L_{\gamma} (\C^n),$ where $\gamma > 0$ such that 
 $u \leq \vphi$ on $\w.$ 
 \end{cor}
 \proof This result follows from the last theorem applied to 
the function $\Tilde{\vphi}_{\w},$ which belongs to 
$\mcal F (\W)$ and satisfies  $\Tilde{\vphi}_{\w} = \vphi$ on 
$\w$ by Proposition 2.2.  
 \fin  

It's possible to get a control on the Monge-Amp\`ere measure of the
subextension in some cases as the following result shows.
Recall that $\ffa (\Om $ is the set of \psh functions $\vphi \in \mcal F (\Om )$
such that $(dd^c \vphi)^n$ puts no mass on pluripolar subsets
of $\W.$
\begin{thm} Let $\W \Subset \C^n$ be a bounded
 hyperconvex domain and $\psi \in \ffa (\Om )$  with
 $\int _{\Om } (dd^c \psi)^n = 1.$ Then there exists
 $u\in \mathcal L (\C^n)$ such that $u\leq \psi$ on $\W$ and
 $\du = {\bf 1_{\W}} (dd^c \psi)^n$ on $\C^n.$
 Here $\du $ is the unique measure with the property that for
 any sequence $v_j \in \mcal L^+ (\C^n) $ decreasing to $u$ we
 have $\ddvj
\to \du$ weakly on $\C^n.$
\end{thm}
\proof Take a hyperconvex domain $\W'$ containing $\ove \W.$
Then by [Ce 3], there exists a \psh function $\fii \in \mcal F^a (\W')$ such
that $\fii \leq \psi$ on $\W$ and
 $ (dd^c \fii)^n = {\bf 1_\W} (dd^c \psi )^n =:  d\mu$ as Borel measures on
 $\W'.$ Then $\mu (U_j) \to 0$ as $j\to\infty$ when $U_j =\{z \in \W' ;
\fii (z) < -j \} .$ Set $d\mu _j = d\mu -d\mu _{| U_j }$ and observe
that $(dd^c \sup\{\fii , - j\})^n \geq  d \mu_j$ on $\W'$
 (see [De], [Ce-Ko]). Then by [Ko 1]
there exist $\fii _j \in \mcal E_0 (\W') $ such that
 $(dd^c \fii _j)^n =d\mu _j.$ Define
$$
\al _j =\frac{1 }{\mu_j (\Om')}$$ and observe that
$\al_j \geq 1.$ By [Ko 2]
there exist $u_j \in  \mcal L^+ (\C^n) $ with $\sup _{\Om' }u_j =0, u_j \leq
\fii _j $ and $\dduj =\al_j d\mu _j .$

Set
$$
v_j = (\sup _{k\geq j} u_k )^{\ast }.
$$
Then $v_j \geq u_j $ and $\psi \geq \fii \geq u:=\lim v_j \in \lp
.$ Observe that for a fixed $j \in \N,$ the sequence 
$$ \tilde v_{j,k} : = \sup \{u_{\ell} ; j \leq \ell \leq k\}, k \geq j$$
is an increasing sequence of \psh functions in $ \mcal L^+ (\C^n) $ which converges a.e. on
$\C^n$ to $v_j.$ Since $(dd^c u_\ell)^n \geq d \mu$ on $\C^n \sm U_j,$ for any $\ell \geq j,$ it
follows from ([De], [Ce-Ko]) that $(dd^c v_{j,k})^n \geq d \mu$ on
$\C^n \setminus U_j.$ By the convergence theorem [Be-Ta 1], it follows that
$$
(dd^c v_j)^n \geq d \mu \text{    on   } \C^n \setminus U_j $$ and for
$M>0$
$$(dd^c \max (v_j ,-M) )^n \geq d\mu \ \ \text{    on   } \
 \C^n \setminus (U_j \cup V_M), \ \ V_M :=\{ u<-M \} .$$
Since by the convergence theorem [Be-Ta 1],
$$
\lim _{j\to\infty } (dd^c \max (v_j ,-M) )^n =(dd^c \max (u ,-M)
)^n $$ we obtain
$$
(dd^c \max (u ,-M) )^n =\lim _{j\to\infty } (dd^c \max (v_j ,-M)
)^n \geq d\mu \text{    on   } \C^n \setminus  V_M .
$$
Therefore
$$\lim _{M\to\infty } (dd^c \max (u ,-M)
)^n \geq d\mu $$
on $\C^n$ and since the integrals of both measures are
equal, the measures themselves are equal. Hence
$$
d\mu _M = (dd^c \max (u ,-M) )^n \to d\mu .
$$
Take a sequence $(w_j)$ of continuous functions in $ \mcal L^+ (\C^n)$ decreasing to $u$. We have
to prove that $\ddwj \to d \mu.$ It is no loss of generality
to assume that $w_{j}>w_{j+1}$ for all $j$. 

Set for $j \in \N,$
$$d\nu _j =  (dd^c w_j )^n,$$
and
 $$d\nu _{j,M} =  (dd^c \max (w_j ,-M) )^n .$$ 
Fix $t>1.$ Then
for $E_j =\{ w_j >u +(t-1)\}\cap \{u\geq -M+1 \}$ we have
\begin{equation}
\int_{E_j}d\mu \to 0.
\label{eq:star}
\end{equation}
(since $E_j$ decrease to $\emptyset .$) 
Since the set $\{ w_{j}<-M\} $ is relatively compact
and the sequence $w_{j}$ is strictly monotone one 
can find $k_{0}$ so big that for $k>k_{0}$, $v_{k
} < w_{j}$ on this set.
Note that if $w_j (z) <
-M$ then $w_j (z) +M >v_{k}(z) +M >
t(v_{k}(z) +M ).$ Hence, by the
comparison principle
$$\aligned
\int _{ \{w_j <-M \} }\ddwj & \leq 
\int _{ \{ t(v_{k}+M) < w_j +M \}
}\ddwj \\ & \leq t^n \int _{ \{ t(v_{k}+M) < w_j +M
 \} } (dd^c v_{k} )^n \leq t^n \int _{ \{u <- M+1 \}\cup
  \ E_j  } (dd^c v_{k} )^n. \endaligned $$
 Then, by
(\ref{eq:star})
$$
\limsup _{j\to\infty }\int _{ \{w_j <-M \} }\ddwj 
\leq \liminf_{k\to\infty } \int _{ \{u <-
M+1 \} } (dd^c v_{k} )^n  =: \ep (M).
$$
From this estimate and the fact that $d\nu _{j,M} =d\nu _j $ on $
\{w_j > -M \}$ we conclude that the total variation
$$
|| \nu _{j,M} -\nu _j ||\leq 2\ep (M/2), \ \ j\geq j(M).
$$

We claim that $\ep (M)\to 0$ as $M\to\infty.$ Indeed, since
$\int_{\C^n} (dd^c v_k)^n = 1 = \mu (\C^n)$
and $(dd^c v_k)^n \geq d \mu$ on $\W \sm U_k
 \supset \W \sm (U_k \cup V_{M - 1}),$ it follows that
 $\int_{U_k \cup V_{M - 1}} (dd^c v_k)^n \leq \mu (U_k \cup V_{M - 1}).$
 Now since $ \mu (U_k
 \cup V_{M - 1}) \leq \mu (U_k) + \mu ( V_{M - 1})$ and
 the measure $\mu$ puts no mass on pluripolar sets, it follow that each of these terms tends
 to $0$ and then so does $\epsilon (M) \leq \mu (V_{M - 1}).$

 Therefore for a test
function $\chi$ we can make the first and the third term on the
right in the formula
$$
\int\chi d(\nu _j -\mu ) = \int\chi d(\nu _j -\nu _{j,M}
)+\int\chi d( \nu _{j,M}-\mu _M )+\int\chi d(\mu _M -\mu )
$$
 arbitrarily small by taking $M$ large enough and 
 $j\geq j(M)$. The middle term
goes to zero as $j\to\infty$ by the convergence theorem. Therefore
the left hand side tends to zero.
\fin
\vskip 0.2 cm
\noindent {\bf Remarks:} 1) Observe that the measure $\du $ was defined
globally. It would be interesing if we can show that it is can be defined
locally.
 It would be also interesting to know if the last theorem is
 true for $\psi
 \in \mathcal F (\Om )$.

\noindent 2) We can define a "canonical" subextention
$$
u =\sup \{ v\in \lp : v\leq\psi \} .
$$
Roughly speaking it should have \MA mass supported on the set
 $ \{ u = \psi \}$.

We will come back to these questions in a subsequent paper.

\vskip 0.3 cm
\noindent {\bf Example:} We give an example of a bounded subharmonic function
$v$ on the unit disc $D \sub \C$ such that $\int_D dd^c v < + \infty$ and
 there is no bounded subharmonic function $u$ on
an open neighbourhood $D'$ of $\ove{D}$ such that $ dd^c u = {\bf 1}_{D} dd^c v$ as
measures on $D'.$

Indeed let $(a_{j})$ be a discreete sequence of points
 in the unit disc $D \sub \C$ wich converges
 to $1.$
 For $j \in \N^{*},$ define $v_{j} (z) :=
 \sup \{ g_{D} (z,a_{j}) , - 1\}$ for $z \in D,$
 where $g_{D} (z,a_{j})$ is the Green function
 of $D$ with pole at $a_{j}.$ Let $(\vep_{j})$ be a
 sequence of positive numbers such that
 $\sum_{j} \vep_{j} = 1.$
 Then $ v := \sum_{j} \vep_{j} v_{j}$ is a bounded
 subharmonic function on $D$ such that
 $ - 1 \leq v \leq 0$ on $D$ and $\int_D dd^c u < + \infty.$

 It is easy to see that $dd^c v_{j} $ converges
 weakly to $\delta_{1}$
as measures on $\C$ since $\int_{D} dd^c v_{j}
 = 1$ for any
 $j \in \N^{*}$ and for $j \in \N$ large enough
 $dd^c v_{j}$ puts no mass outside any arbitrary
 neighbourhood of $1.$
 Then it follows that 
 $\limsup_{j \to + \infty} \int_{\ove{D}}
 \log \vert 1 - z\vert  dd^c v_{j} =
  - \infty.$
 Therefore, taking a subsequence if necessary, we can assume
 that the poles and the weights are choosen so that
 $\sum_{j} \vep_j \int_{\ove{D}}
 \log \vert 1 - z\vert  dd^c v_{j} = - \infty.$

 Now if $u $ is a subharmonic function on a disc
 $D'$ containing $\ove{D},$ such that $dd^c u =
 {\bf 1}_{D} dd^c v$ on $D'$ then by the Riesz decomposition we have 
$ u (1) = c +
 \int_{\ove D} \log \vert 1 - z \vert  dd^c u = c +
 \sum_{j} \vep_j \int_{\ove{D}}
  \log \vert 1 - z \vert  dd^c v_{j} = - \infty $ by
 construction.

\vskip 0.3 cm
\noindent {\bf Remarks:}
 As observed  by El Mir ([El]),  there exists a \psh function 
 $u $ on some open subset $\W \sub \C^2$ which  has no 
 subextension to any  larger domain.  The main obstruction
  is  the  fact that the polar set of $u$ in $\W$ may contain a non 
 trivial analytic set, which does not extend as an analytic set in a 
 larger domain. This  analytic stucture comes from the fact that
  superlevel sets of Lelong numbers of $\vphi$ defined by 
 $$A (\vphi;c) := \{a \in \W \ ; \  \nu (\vphi,a) \geq c\}, c > 0,$$
 are analytic sets by Siu's theorem (cf. [Siu]).
  Indeed, let $u$ be a \psh function on  some  open subset 
$\omega  \sub \W.$ Assume that there exists a function 
$U \in PSH (\W)$ such that $U \leq  u $ on $\omega.$ Then  
 $A (u;c)  \subset  A (U;c) \cap \omega.$ 
Hence the subextension problem is closely related to the propagation of singularities of
plurisubharmonic functions.
Observe that for all the functions which were considered in our theorems the sets $ A
(\vphi;c)$ are finite so that obviously there is no analytic obstruction to subextension.

\bigskip
Urban CEGRELL \\
Department of Mathematics \\
University of Umea \\
S-90187 Umea, Sweeden \\
and \\
Mid Sweden University\\
Department of Mathematics\\
S-85 170 Sundsvall, Sweeden

\vskip 0.3 cm
\noindent Slawomir KOLODZIEJ \\
Jagiellonian University \\
Institute of Mathematics \\
Reymonta 4, 30-059 Krak\'ow, Poland
\vskip 0.3 cm

\noindent Ahmed ZERIAHI \\
Universit\'e Paul Sabatier \\
Institut de Math\'ematiques \\
118 Route de Narbonne \\
31062 Toulouse cedex, France

 \end{document}